\documentclass[12pt]{amsart}
\usepackage{amsmath,amsfonts,amssymb,amsthm}
\usepackage{pinlabel}
\usepackage{graphicx}
\usepackage{mathtools}
\usepackage[all]{xy}
\usepackage{hyperref}
\usepackage[usenames,dvipsnames]{color}
\usepackage{color}
\xyoption{dvips}

\newtheorem{theorem}{Theorem}[section]

\newtheorem{proposition}[theorem]{Proposition}

\theoremstyle{remark}
\newtheorem{remark}[theorem]{Remark}

\theoremstyle{example}
\theoremstyle{definition}

\newcommand{\R}{\mathbb{R}}
\newcommand{\Z}{\mathbb{Z}}

\newcommand{\A}{\mathcal A}

\newcommand{\ve}{\varepsilon}

\numberwithin{equation}{section}

\begin{document}

\title[On torsion in linearized Legendrian contact homology]{On torsion in linearized Legendrian contact homology}

\author{Roman Golovko}

\begin{abstract}
In this short note we discuss certain examples of Legendrian submanifolds, whose linearized Legendrian contact (co)homology groups over integers have non-vanishing algebraic torsion. More precisely, for a given arbitrary finitely generated abelian group $G$ and a positive integer $n\geq 3$, $n\neq 4$, we construct  examples of Legendrian submanifolds of the standard contact vector space $\R^{2n+1}$, whose $n-1$-th linearized Legendrian contact (co)homology over $\mathbb Z$ computed with respect to a certain augmentation is isomorphic to $G$.
\end{abstract}

\address{Faculty of Mathematics and Physics, Charles University, Sokolovsk\'{a} 49/83, 18675 Praha 8, Czech Republic}
\email{golovko@karlin.mff.cuni.cz}
\date{\today}
\thanks{}


\maketitle

\section{Introduction and main result}
The Legendrian contact homology of a closed Legendrian submanifold $\Lambda$ of the standard contact vector space $(\R^{2n+1},\xi_{st}=\ker \alpha_{st})$, where $\alpha_{st}=dz-ydx$, is a modern Legendrian invariant defined by Eliashberg--Givental--Hofer \cite{EGHSFT} and Chekanov \cite{ChekanovDGAL}, and developed by Ekholm--Etnyre--Sullivan \cite{EkhEtnSul2005}.
It is a homology of the Legendrian contact homology (LCH) differential graded algebra (often called the Chekanov--Eliashberg differential graded algebra). Chekanov--Eliashberg  DGA is a unital noncommutative differential graded algebra freely generated by the (generically finite) set of integral curves of the Reeb vector field $\partial_z$ that start and end on $\Lambda$. These integral curves are called  Reeb chords of $\Lambda$.
Legendrian contact homology is often defined over $\Z_2$, but if $\Lambda$ is spin it can be also defined over other fields, over $\Z$ \cite{OrientLegContHom, KarlssonLCHor} and even more general coefficient rings whose structure envolves certain topological information about $\Lambda$ such as $\Z_2[H_1(\Lambda;\Z)]$ or $\Z[H_1(\Lambda;\Z)]$ \cite{EkhEtnSul2005, KarlssonLCHor}. 

The Legendrian contact homology DGA is not finite rank, even in fixed degree; the same holds in homology: the
graded pieces of the Legendrian contact homology are often infinite dimensional and difficult to compute.
In order to deal with this issue Chekanov \cite{ChekanovDGAL}  proposed to use an augmentation of
the DGA to produce a generically finite-dimensional linear complex, whose homology is called linearized Legendrian contact homology.

Given an exact Lagrangian filling $L$ of $\Lambda$ in the symplectization $(\R\times \R^{2n+1}, d(e^t \alpha_{st}))$ with vanishing Maslov number, it  induces
the augmentation of the Chekanov--Eliashberg algebra, i.e. a unital DGA homomorphism $\varepsilon: \mathcal A(\Lambda)\to (\Z_2,0)$, see \cite{EkhomHonaKalmancobordisms}. 
If besides that $L$ is equipped with a spin structure extending the given spin structure on $\Lambda$, then one also has an augmentation $\varepsilon: \mathcal A(\Lambda)\to (\Z,0)$, see \cite{EkhomHonaKalmancobordisms, KarlssonOrc}.

Most of the computations of linearized Legendrian contact homology groups have been done for the Chekanov-Eliashberg algebras with $\Z_2$-coefficients. One can ask whether for integral coefficients 
one can get some elements of finite order in the linearized Legendrian contact (co)homology, i.e. if one can get a non-trivial algebraic torsion in linearized Legendrian contact (co)homology. Note that certain examples of Legendrian submanifolds of $J^1(T^ n)$ with torsion elements in linearization appear in the work of Ekholm--Etnyre--Sullivan \cite[Theorem 5.4]{OrientLegContHom}. Besides that, in the work of Ng on knot contact homology \cite{NgKBIFCH1, NgKBIFCH2} there are  examples of Legendrians in $J^1(S^2)$ with certain torsion elements in linearization. The question that we would like to answer is about the type of algebraic torsion that could appear for  Legendrian submanifolds of the standard contact vector space.  More precisely, we would like to answer the question of whether an arbitrary 
finitely generated abelian group can be realized as a linearized Legendrian contact (co)homology of some Legendrian submanifold of the standard contact vector space.  

We provide the following answer to this question in high dimensions:
\begin{theorem}
\label{topresnongaug}
Given a finitely generated abelian group $G$ and an integer $n\geq 3$, $n\neq 4$.
There is a Legendrian submanifold $\Lambda$ in $\mathbb R^{2n+1}$ of Maslov number $0$ such that the Chekanov-Eliashberg algebra of $\Lambda$  admits an augmentation 
$\varepsilon: \A(\Lambda)\to (\Z, 0)$ with $LCH_{\varepsilon}^{n-1}(\Lambda;\Z)\simeq G$.
\end{theorem}

\section{Proof of Theorem \ref{topresnongaug}}

We start with the following construction of a spin manifold, whose first homology with $\Z$-coefficients is isomorphic to $G$.
\subsection{Construction of a spin manifold}
\label{spincontrolledh1}
Given a finitely generated abelian group $G$, we can write it as 
$$G\simeq \Z^k\times \Z_{p_1}\times\dots 
\times  \Z_{p_s},$$
where $k$ is a non-negative integer and $p_1,\dots, p_s$ are non-negative  integers that are powers of some not necessarily distinct prime numbers. 
Then we construct a closed, oriented, connected $3$-manifold $N$ such that $H_1(N;\mathbb Z)\simeq G$. We get $N$ from the $k$ distinct copies of $S^1\times S^2$ and the collection of lens spaces $L(p_1,q_1),\dots,L(p_s,q_s)$. More precisely, since we know that
\begin{align}
H_1(S^1\times S^2;\mathbb Z)\simeq \Z\quad \mbox{and}\quad H_1(L(p_i,q_i);\Z)\simeq \Z_{p_i},
\end{align}
we define $$N= \#^k (S^1\times S^2)\# L(p_1,q_1)\# \dots \#L(p_s,q_s).$$
Since we know that all closed, orientable manifolds of dimension $3$ are spin, we can say that $N$ is spin. 
We observe that 
\begin{align*}
&H_1(N;\Z)=H_1(\#^k (S^1\times S^2) \# L(p_1,q_1)\# \dots \#L(p_s,q_s);\Z)\\  &\simeq H_1(S^1\times S^2;\Z)^k\times   H_1(L(p_1,q_1);\Z)\times \dots \times H_1(L(p_s,q_s);\Z)
\\ &\simeq\Z^k\times \Z_{p_1}\times\dots 
\times  \Z_{p_s}\simeq G.
\end{align*}
Then since $S^n$ is spin for all $n\geq 0$ and the product of spin manifolds is a spin manifold, we can say that $N\times S^{n-3}$ is a spin manifold. 
Finally, when $n\neq 4$, K\"{u}nneth formula implies that 
\begin{align}
\label{compihom}
H_{1}(N\times S^{n-3}; \Z)\simeq H_{1}(N; \Z)\simeq G.
\end{align}
From now on we will denote $N\times S^{n-3}$ by $\Lambda'$.
\subsection{Construction of a Legendrian submanifold and its filling}
Now we recall the following statement, which can be seen as an implication of the Murphy's h-principle, for the justification see \cite[Appendix A, Remark A.5]{LLEIHDCM} and \cite[Theorem  7.9, Lemma 7.17]{CielebakEliashberg2012}.
\begin{proposition}
\label{PorpGeoRealization}
Every closed spin $n$-dimensional manifold admits a Legendrian embedding  $i:\Lambda \to (\R^{2n+1}, \xi_{st})$  with a vanishing Maslov number.  
\end{proposition}
\begin{remark}
 Str\"angberg in \cite{LegendrianAprroximation} described the $C^0$-approximation procedure for $2$-dimensional submanifolds in the standard contact $\R^5$ by Legendrian submanifolds generalizing the standard approximation result for closed Legendrians in the standard contact $\R^3$. The construction of  Str\"angberg as mentioned in \cite{LegendrianAprroximation} admits an extention to certain high dimensions,  that will also potentially lead to Proposition \ref{PorpGeoRealization}.
\end{remark}
We apply Proposition \ref{PorpGeoRealization} to $\Lambda'$ from Section \ref{spincontrolledh1} and get a Legendrian embedding  $i:\Lambda' \to (\R^{2n+1}, \xi_{st})$  with a vanishing Maslov number, this Legendrian embedding will still be denoted by $\Lambda'$. We push  $\Lambda'$ slightly in $\partial_z$-direction (Reeb direction) and get a $2$-copy of $\Lambda'$ that we call $\Lambda$. Next we apply the technique for constructing Lagrangian cobordisms
from the Reeb flow discussed by Mohnke in \cite{Mohnketorusfilliingpair}. In particular, the technique of Mohnke in the deformed form written by Dimitrolglou Rizell in \cite[Lemma 2.5]{DimitroglouRizell2023} applied to $\Lambda'$ and a $U$-shaped curve implies that  $\Lambda$ admits an exact Lagrangian filling of Maslov number $0$ diffeomorphic to $\R\times \Lambda'$.
Following \cite{EkhomHonaKalmancobordisms,KarlssonOrc}, we observe that this spin Maslov number $0$ exact Lagrangian filling induces the augmentation $\varepsilon$ to $\mathbb Z$. Then we would like to study linearized Legendrian contact cohomology of $\mathcal A(\Lambda)$,  linearized with respect to $\varepsilon$.

\subsection{Application of the isomorphism of Seidel-Ekholm-Dimitroglou Rizell}

Recall that there is an isomorphism due to Seidel  that has been described
by Ekholm in \cite{Seidelsisowfc} and completely proven by Dimitroglou Rizell in \cite{DimitroglouRizellLiftingPSH}.
\begin{theorem}[Seidel--Ekholm--Dimitroglou Rizell]
\label{SEDRI}
Let $\Lambda$ be a Legendrian submanifold of Maslov number $0$ of $\R^{2n+1}_{st}$, which admits an exact Lagrangian filling $L$ of Maslov number $0$. Then
\begin{align}
LCH_{\ve}^{i}(\Lambda; \Z_{2})\simeq H_{n-i}(L_{\Lambda}; \Z_{2})
\end{align}
where $\varepsilon$ is the augmentation induced by $L$.
\end{theorem}

The homology and cohomology groups in
the above result are defined over $\Z_2$. 
Recall that after the proof of Dimitroglou Rizell \cite{DimitroglouRizellLiftingPSH} signs of the cobordism maps between Chekanov-Eliashberg algebras 
have been studied by Karlsson in \cite[Theorem 2.5, Theorem 2.6]{KarlssonOrc}. In addition, the work of Ekholm--Lekili \cite[Section 1.2]{EkholmLekiliduality} implies the enhancement of Theorem \ref{SEDRI}, which compares not just the corresponding homology and cohomology groups, but the
corresponding $A_{\infty}$-structures, which holds with signs and works over an arbitrary field. Therefore, we can say that the isomorphism of Seidel--Ekholm--Dimitroglou Rizell holds over $\Z$.

Then we take the linearized Legendrian contact cohomology of $\Lambda$ and apply the isomorphism of Seidel-Ekholm--Dimitroglou Rizell over $\Z$ and get that
\begin{align}
\label{Seidelsisofor2copy}
LCH^{n-1}_{\varepsilon}(\Lambda; \mathbb Z) \simeq H_{1}(\mathbb R\times \Lambda'; \Z)\simeq H_{1}(\Lambda';\Z).
\end{align}
Then combining Formula (\ref{Seidelsisofor2copy}) with  Formula (\ref{compihom}) we get that
\begin{align}
\label{thefinform}
LCH^{n-1}_{\varepsilon}(\Lambda; \mathbb Z)\simeq G.
\end{align}

\begin{remark}
We can also take the $1$-jet space $J^1(\Lambda')\simeq \R\times T^{\ast} \Lambda'$ equipped with the standard contact $1$-form $\alpha=dz+\theta$, where 
$\theta$ is a primitive of the standard symplectic form on $T^{\ast} \Lambda'$. Then we take $0_{\Lambda'}$, and denote
by  $\Lambda$ the $2$-copy of $0_{\Lambda'}$. In this case one does not need to apply the isomorphism of  Seidel--Ekholm--Dimitroglou Rizell, but 
one can simply rely on the analysis of pseudoholomorphic curves in the duality paper of Ekholm--Etnyre--Sabloff  \cite{EESDuality} to get Formula \ref{thefinform}), and hence the analogue of Theorem \ref{topresnongaug} will hold.
\end{remark}

\begin{remark}
Note that the main question of this paper can be naturally extended to the question on geography of (bi)linearized Legendrian contact (co)homology groups over $\Z$. In particular, one can ask whether an arbitrary graded collection of finitely generated abelian groups can be realized as a graded collection of (bi)linearized Legendrian contact (co)homology groups of a given Legendrian under the condition that these groups naturally fit into the duality long exact sequence \cite{EESDuality}. This type of question for bilinearized Legendrian contact (co)homology over $\Z_2$ has been studied by Bourgeois--Galant \cite{BourgeoisGalant19}, and for generating families (co)homology by Bourgeois--Sabloff--Traynor \cite{BourgeoisSabloffTraynor2015}.
\end{remark}

\begin{remark}
Note that the examples  we construct have ``complicated'' topology that allows us to realize an arbitrary finitely generated abelian group as a linearized Legendrian contact (co)homology group. This type of topology exists only in high dimensions. Our construction is inapplicable for Legendrian knots and links. 
On the other hand, one could ask whether an arbitrary finitely generated abelian group can be realized as a linearized Legendrian contact (co)homology of a Legendrian with ``simple'' topology, for example Legendrian knots and links. This question is due to Bourgeois \cite{BourgeoisTorsion} and it remains wide open. 
\end{remark}

\section*{Acknowledgements}
The author is grateful to Russell Avdek, Georgios Dimitroglou Rizell and  Filip Strako\v{s} for the very helpful discussions. 
Besides that, the author would like to thank the anonymous referee for suggesting quite a few useful improvements.
The author is supported by the GA\v{C}R EXPRO Grant 19-28628X.

\end{document}